\documentclass[aap,MSNbibl,citesort,seceqn,dvips]{arximspdf}

\doi{10.1214/12-AAP856} 
\volume{23}
\issue{2}
\pubyear{2013}
\firstpage{811}
\lastpage{833}

\makeatletter
\newcommand{\implies}{\Longrightarrow}
\newcommand{\Yv}{\mathbf{Y}}
\newcommand{\Xv}{\mathbf{X}}
\newcommand{\ev}{\mathbf{e}}
\newcommand{\av}{\mathbf{a}}
\newcommand{\pv}{\mathbf{p}}
\newcommand{\qv}{\mathbf{q}}
\newcommand{\zv}{\mathbf{z}}
\newcommand{\uI}{\mathbf{I}}
\newcommand{\uM}{\mathbf{M}}
\newcommand{\uP}{\mathbf{P}}
\newcommand{\ualpha}{\bolds{\alpha}}

\newtheorem{thmm}{Theorem}[section]
 \newtheorem{prop}{Proposition}[section]
 \newproclaim{eg}{Example}[subsection]
\newproclaim{rem}{Remark}
\newcommand{\E}{\mathrm{E}}
\newcommand{\Po}{\mathrm{P}}
\makeatother

\begin{document}
\begin{frontmatter}

\title{Convergence analysis of some multivariate Markov chains
using stochastic monotonicity}
\runtitle{Convergence analysis of Markov chains}

\begin{aug}
\author[A]{\fnms{Kshitij} \snm{Khare}\corref{}\ead[label=e1]{kdkhare@stat.ufl.edu}}
\and
\author[B]{\fnms{Nabanita} \snm{Mukherjee}\ead[label=e2]{mukherjeen@email.chop.edu}}
\runauthor{K. Khare and N. Mukherjee}
\affiliation{University of Florida and Center for Outcome Research}
\address[A]{Department of Statistics\\
University of Florida\\
Gainesville, Florida 32611\\
USA\\
\printead{e1}}

\address[B]{Center for Outcome Research\\
The Children's Hospital of Philadelphia\\
3535 Market Street, Suite 1005\\
Philadelphia, Pennsylvania 19104\\
USA\\
\printead{e2}}

\end{aug}

\received{\smonth{5} \syear{2011}}
\revised{\smonth{2} \syear{2012}}

%
\begin{abstract}
We provide a nonasymptotic analysis of convergence to stationarity for
a collection of
Markov chains on multivariate state spaces, from arbitrary starting
points, thereby
generalizing results in [Khare and Zhou \textit{Ann. Appl. Probab.}
\textbf{19} (2009) 737--777].
Our examples include the multi-allele Moran
model in population genetics and its variants in community ecology, a
generalized
Ehrenfest urn model and variants of the P\'{o}lya urn model. It is
shown that all these
Markov chains are stochastically monotone with respect to an
appropriate partial ordering.
Then, using a generalization of the results in
[Diaconis, Khare and Saloff-Coste \textit{Sankhya} \textbf{72} (2010)
45--76] and
[Wilson \textit{Ann. Appl. Probab.} \textbf{14} (2004) 274--325] (for
univariate totally ordered spaces) to multivariate partially
ordered spaces, we obtain explicit nonasymptotic bounds for the
distance to stationarity
from arbitrary starting points. In previous literature, bounds, if any,
were available
only from special starting points. The analysis also works for
nonreversible Markov
chains, and allows us to analyze cases of the multi-allele Moran model
not considered
in [Khare and Zhou \textit{Ann. Appl. Probab.} \textbf{19} (2009) 737--777].
\end{abstract}

%
\begin{keyword}[class=AMS]
\kwd[Primary ]{60J10}
\kwd[; secondary ]{60J22}
\kwd{47H05}.
\end{keyword}

\begin{keyword}
\kwd{Convergence analysis}
\kwd{Markov chains}
\kwd{stochastic monotonicity}
\kwd{partial ordering}.
\end{keyword}

\end{frontmatter}
%

\section{Introduction}\label{sec1}

The theory of Markov chains plays a prominent role in the fields of
statistics and applied probability. Markov chains
have a wide range of applications in numerous areas from particle
transport through finite state machines to the
theory of gene expression. Some important applications include modeling
scientific phenomena in population genetics,
statistical physics and image processing. Another important use is
simulating from an intractable probability
distribution. It is a well-known fact that, under mild conditions
discussed in~\cite{athdosethu}, a Markov chain
converges to its stationary distribution. In the applications mentioned
above, often it is useful to know exactly
how long to run the Markov chain until it reaches sufficiently close to
the stationary distribution. Answering this
question as accurately as possible, is what obtaining a
``nonasymptotic convergence analysis'' of Markov chains is all
about. The applied probability community has made significant strides
in this area in the past three decades. Despite
this progress, answering this question still remains a challenging task
for various standard Markov chains arising in
applied probability and statistics. There are various examples where
currently available state of the art techniques
can give upper bounds that are substantially larger than the correct
answer, often by orders of magnitude.

In the current paper, we provide a nonasymptotic analysis of
convergence to stationarity for a collection of Markov
chains in population genetics. The analysis is based on a
generalization of the monotone coupling argument to
multivariate state spaces. These Markov chains appear as standard
models in population genetics and ecology and
include the multi-allele Moran process in population genetics and its
variants in community ecology, a generalized
Ehrenfest urn model and the P\'{o}lya urn process. These Markov chains
were analyzed in~\cite{Khare09}, and the
authors provide an exact convergence analysis in terms of the
``chi-square distance'' by using spectral techniques.
But their analysis is somewhat incomplete because it works \textit{only
for some natural selected starting points}.
Stochastic monotonicity of a Markov chain, along with the knowledge of
a monotone eigenfunction (see~\cite{Diaconis06Coupling} and~\cite{wncoupling}), can be used to obtain
a nonasymptotic convergence analysis from
an arbitrary starting point. Existing results in \cite
{Diaconis06Coupling} and~\cite{wncoupling} require total
ordering of the state space, which generally works in the case of
univariate state spaces. In multivariate state
spaces, however, there often exists a natural partial ordering. We
prove that the Markov chains being considered in
this paper are stochastically monotone with respect to an appropriate
partial ordering; see Theorems~\ref{them2},
\ref{them3},~\ref{them4}. \textit{But stochastic monotonicity of a Markov
chain with respect to the partial ordering,
even with the knowledge of a monotone eigenfunction, is not enough to
get desired convergence bounds.} However, an
additional condition, satisfied by all the Markov chains under
consideration in this paper, enables us to obtain
useful convergence bounds; see Theorem~\ref{them1}. Another limitation
of the spectral techniques used in
\cite{Khare09} is that they require reversibility of the Markov chain
under consideration. The coupling argument
presented in this paper also works for nonreversible Markov chains.
Using this, for example, we are able to obtain
explicit convergence bounds for generalizations of the standard
multi-allele Moran model which are nonreversible.

Another important issue to understand is that out of the three classes
of examples considered in this paper, the
stationary distribution and the second largest eigenvalue of the Markov
chains corresponding to the generalized
Ehrenfest urn models and the P{\'o}lya urn models are known (\textit{the
stationary distribution is unknown for the
general multi-allele Moran model}). Hence, for these two models, from a
general starting point $\mathbf{x}$, one
could potentially consider the crude upper bound $\frac{\lambda^n}{2
\sqrt{\pi(\mathbf{x})}}$ for the total variation
distance from stationarity after $n$ steps. Here $\pi(\mathbf{x})$ denotes
the mass put by the stationary distribution
at $\mathbf{x}$, and $\lambda$ denotes the second largest eigenvalue. However,
the upper bounds derived in this paper
mostly provide a significant improvement over the crude upper bound.
See the remarks in Section~\ref{polyabound}
and Section~\ref{ehrenbound}.

Here is an example of our results. The Unified Neutral Theory of
Biodiversity and Biogeography (UNTB) is an important
theory proposed by ecologist Stephen Hubbell in his monograph \cite
{Hubbell01} which is used in the study of
diversity and species abundances in ecological communities. There are
two levels in Hubbell's theory, a~metacommunity
and a local community.

We concentrate here on the evolution of the local community. The local
community has constant population size $N$ with
$d$ different species. At each step, one individual is randomly chosen
to die and is replaced by a new individual.
With probability $m$, the new individual is chosen randomly from the
metacommunity, which has proportion $p_i$ of
species $i$ ($i=1,2,\ldots,d$). With probability $1-m$, the new
individual is randomly chosen from the
remaining $N-1$ individuals in the local community. This process is a
variant of the so-called multi-allele Moran
model in population genetics~\cite{Ewens04}. The metacommunity evolves
at a much larger time scale and is assumed to
be fixed during the evolution of the local community.

A very important issue of both practical and theoretical interests is
to determine how soon a local community reaches
equilibrium (see McGill~\cite{McGill03}). Let $K(\cdot,\cdot)$ be the
transition density of our local community Markov
chain with state space~$\mathcal{X}$ and stationary density $\pi$. Let
$\mathbf{x}\in\mathcal{X}$ be the initial state of
the Markov chain. We are interested in answering the following
question. For arbitrary $\varepsilon>0$, how many steps,
$n$, are needed so that the total variation distance between the
density of the Markov chain after $n$
steps and the stationary density is less than $\varepsilon$? More
precisely, we want to find $n$ such that
\[
\| K_\mathbf{x}^n - \pi\|_{\mathrm{TV}} = \frac{1}{2} \sum
_{\mathcal{X}}
|K^n (\mathbf{x}, \mathbf{x}^\prime) - \pi(\mathbf{x}^\prime)|
\leq
\varepsilon,
\]
where $K^n_{\mathbf{x}}$ denotes the density of the chain started at
state $\mathbf{x}
$ after $n$ steps.\footnote{For ease of
exposition, if $f$ and $g$ are densities with respect to the counting
measure on a finite state space $\mathcal{X}$,
$\|f - g\|_{\mathrm{TV}}$ will denote the total variation distance
between the probability measures
corresponding to $f$ and $g$.}

Khare and Zhou~\cite{Khare09} provide an exact answer to this question
in terms of the ``chi-square distance''
by using spectral techniques, \textit{when all individuals belong to the
same species to begin with}. So providing any
nonasymptotic convergence bounds from an arbitrary starting point was
still unresolved. Convergence bounds for the
general local community Markov chain are provided in Section~\ref{secmoran},
with an arbitrary starting point. Note
that the upper and lower bounds obtained are not exactly matching, but
they are within a reasonable range of each
other. Considering the fact that no useful analysis was available from
an arbitrary starting point, the bounds
provided are definitely a significant step forward.

As an illustration, note that under suitable parametrization (see \cite
{Khare09}), the local community process by
Hubbell is the same as the P\'{o}lya down--up model; see Section \ref
{polya1}. Suppose that the local community
has population size \mbox{$N=100$} with \mbox{$d=5$} species. With probability
$m=0.9$, the new individual is chosen randomly
from the meta-community with uniform species frequencies $\mathbf{p} =
(0.2, 0.2, 0.2, 0.2, 0.2)$. Let
$\Xv=(X_1,\ldots,X_d)$ be any (random) count vector of the local
community, where $X_i$ is the count of
individuals of species $i$. From Section~\ref{polya1}, for a starting
state $\mathbf{x}=(0,10,0,10,80)$, the bounds on
the total variation distance are obtained as
%
\begin{equation}\label{eqpolya}
0.375 \biggl( 1 - \frac{1}{111}\biggr)^{n} \leq\| K_\mathbf{x}^n
-\pi
\|_{\mathrm{TV}} \leq100 \biggl( 1 -
\frac{1}{111}\biggr)^{n}.
\end{equation}
For $\varepsilon=0.01$, (\ref{eqpolya}) tells us that at least $401$
steps are necessary and at most $1018$ steps
are sufficient for the total variation distance to be less then $0.01$.
The crude upper bound for total
variation distance is $(2.2186 \times10^{19}) ( 1 - 1/111
)^{n}$ which gives $5432$ steps are
sufficient for the total variation distance to be less then $0.01$.

The paper is organized in the following way. In Section~\ref{sec2}, we
provide the necessary
background for stochastic monotonicity, and then proceed to prove
Theorem~\ref{them1}, which
generalizes the results in~\cite{Diaconis06Coupling} and \cite
{wncoupling} to multivariate
partially ordered finite state spaces to obtain convergence bounds,
under appropriate
monotonicity assumptions. In Section~\ref{sec3}, three classes of
Markov chains:
multi-allele Moran model, generalized Ehrenfest urn model and
generalized P{\'o}lya urn
model are considered. Each of these Markov chains is shown to be
stochastically monotone
with respect to an appropriate partial ordering, and also shown to
satisfy the other
assumptions in Theorem~\ref{them1}. All these are combined to provide
nonasymptotic
convergence bounds for these classes of Markov chains from arbitrary
starting points. We
conclude the paper with a short discussion in Section~\ref{sec4}.

\section{Monotone Markov chains}\label{sec2}

\subsection{Background}\label{sec2.1}

Let $\mathcal{X}$ be a finite state space with total ordering $\leq$.
Let $K(\cdot,\cdot)$ be a Markov kernel on
$\mathcal{X}$. We say $K$ is stochastically monotone if for all $x
\in
\mathcal{X}$ and $x^\prime\in
\mathcal{X}$ with $x \leq x^\prime$,
\[
\sum_{y\leq y^\prime}K(x, y) \geq\sum_{y\leq y^\prime} K(x^{\prime},
y) \qquad  \mbox{for all }   y^\prime\in
\mathcal{X}.
\]
Monotone Markov chains have been thoroughly studied and applied. See
Lund and Tweedie~\cite{Lund96},
Stoyan~\cite{Stoyan83} and the references therein. They are currently
popular because of ``coupling from the
past.'' See David Wilson's website on perfect sampling,
\texttt{%
\href{http://research.microsoft.com/en-us/um/people/dbwilson/exact}{http://research.microsoft.com/en-us/um/}
\href{http://research.microsoft.com/en-us/um/people/dbwilson/exact}{people/dbwilson/exact}}, for
extensive references on this subject.

Alternatively, if the state space $\mathcal{X}$ of a Markov chain is
totally ordered (e.g., a~subset of
$\mathbb{Z}$ and $\mathbb{R}$), then the Markov chain with
corresponding transition operator $K$ is stochastically
monotone if for every monotone function $f\dvtx  \mathcal{X}\rightarrow
\mathbb{R}$, the function $Kf$ is also
monotone. There is a standard coupling technique available for monotone
Markov chains on totally ordered spaces.
Wilson~\cite{wncoupling} uses this coupling technique in the presence
of an explicit eigenfunction to provide
general convergence bounds for stochastically monotone Markov chains on
totally ordered finite state spaces.
Diaconis, Khare and Saloff-Coste~\cite{Diaconis06Coupling} provide
extensions for general state spaces and use
these results to analyze certain two-component Gibbs samplers.

However, for multivariate state spaces, there is often no natural total
ordering, but there exists a natural
partial ordering. For example, if $\mathcal{X}$ consists of
$d$-dimensional vectors, then entry-wise domination
gives rise to a standard partial ordering. A Markov chain with
corresponding transition operator $K$ is monotone
with respect to a partial ordering, if whenever $f\dvtx  \mathcal
{X}\rightarrow\mathbb{R}$ is monotone with respect
to the partial ordering, $Kf$ is monotone with respect to the partial
ordering. See
Fill and Machida~\cite{Fill01}, Beskos and Roberts~\cite{bsksrbttgd},
Roberts and Rosenthal~\cite{rbrtsrsntl} and
the references therein for varied applications. The literature on
perfect sampling mainly consists of various
techniques for simulating from specific distributions on partially
ordered spaces with a unique minimal and
maximal element; see Propp and Wilson~\cite{Propp98}. \textit{Note that,
unlike perfect sampling, our focus is
to analyze given Markov chains corresponding to specific models, and
not to devise Markov chains to simulate from
a specified distribution}.

The theorem listed below generalizes earlier results in Wilson \cite
{wncoupling} and
Diaconis, Khare and Saloff-Coste~\cite{Diaconis06Coupling} (for
univariate totally ordered
spaces) to multivariate partially ordered spaces in order to obtain
nonasymptotic
convergence results.

\subsection{Convergence of monotone Markov chains: General result}\label{sec2.2}

\begin{thmm} \label{them1}
Let $K$ be the transition density of a Markov chain on a finite state
space $\mathcal{X}$ equipped
with a partial ordering, $\preceq$. Suppose that $K$ has a stationary
distribution with density
$\pi$, and the following conditions are satisfied:
\begin{longlist}[(a)]
\item[(a)] $K$ is monotone with respect to the partial ordering,
$\preceq$.
\item[(b)] \label{pwd} (Pair-wise dominance property) For an arbitrary
$\mathbf{x}$ and $\mathbf{y}$ in
$\mathcal{X}$, there exists $\zv(\mathbf{x},\mathbf{y})$ (depends
possibly on $\mathbf{x}$
and $\mathbf{y}$) such that
$\zv$ either dominates $\mathbf{x}$ and $\mathbf{y}$ or is dominated by
both $\mathbf{x}$
and $\mathbf{y}$ with respect
to $\preceq$.\vadjust{\goodbreak}
\item[(c)] $\lambda\in(0,1)$ is an eigenvalue of $K$ with strictly
monotone eigenfunction
$f$ such that
\[
c_1 = \inf_{\mathbf{x}^* \preceq\mathbf{y}^*, \mathbf{x}^* \neq
\mathbf{y}^*}
\{ f(\mathbf{y}^*) -
f(\mathbf{x}^*) \vert
  \mathbf{x}^*, \mathbf{y}^* \in\mathcal{X} \} > 0, \qquad  c_2 =
\sup
_{\mathbf{x}\in
\mathcal{X}} |f(\mathbf{x})| >
0.
\]
\end{longlist}
Then for any starting state $\mathbf{x}$,
\[
\frac{\lambda^{n}}{2c_2} \vert f(\mathbf{x}) \vert\leq\Vert
K^{n}_{\mathbf{x}} -
\pi\Vert_{\mathrm{TV}} \leq
\frac{\lambda^{n}}{c_1} \E\vert f(\Yv) + f(\mathbf{x}) - 2f(\zv
(\mathbf{x},\Yv))
\vert,
\]
where $\Yv\sim\pi$.
\end{thmm}

\begin{pf}
Let $\mathbf{x}^* \in\mathcal{X}$ and $\mathbf{y}^* \in\mathcal{X}$
satisfy $\mathbf{x}
^*\preceq\mathbf{y}^*$. It is
well known that if a probability distribution $\mu$ on $\mathcal{X}$ is
stochastically
dominated by another probability distribution $\nu$ on $\mathcal{X}$,
that is, $\int f \,d\mu
\leq\int f\, d\nu$ for every monotone function $f$, then we can
construct random variables
$\Xv$ and $\Yv$ such that $\Xv\sim\mu,  \Yv\sim\nu$ and $\Xv
\preceq\Yv$; see for
example~\cite{stcinokko}. Since $K$ is monotone with respect to the
partial ordering,
$\preceq$, by repeated application of this result, we can construct two
coupled Markov
chains, $\{\Xv_n\}_{n\geq0}$ and $\{\Yv_n\}_{n\geq0}$ such that
$\Xv
_0=\mathbf{x}^*, \Yv_0=\mathbf{y}^*$
and $\Xv_n\preceq\Yv_n$ for every $n\geq1$. Further, if $\Xv_{n_0} =
\Yv_{n_0}$, then
$\Xv_{n} = \Yv_{n}$ for all $n\geq n_0$.

It follows that for any $n\geq1$,
\begin{eqnarray*}
\Vert K^{n}_{\mathbf{x}^*}-K^{n}_{\mathbf{y}^*}\Vert_{\mathrm{TV}}
&\leq&
\Po(\Xv
_n\neq\Yv_n\vert\Xv_0=\mathbf{x}^*, \Yv_0=\mathbf{y}^*)\\
&\leq& \E\biggl\{\frac{f(\Yv_n)-f( \Xv_n)}{c_1}\Big\vert\Xv
_0=\mathbf{x}^*, \Yv
_0=\mathbf{y}^* \biggr\}.
\end{eqnarray*}
The previous inequality uses $\Xv_{n} \preceq\Yv_{n}$, the strict
monotonicity of $f$ and the
hypothesis that $f(\mathbf{y})-f(\mathbf{x})\geq c_1$ if $\mathbf
{x}\preceq
\mathbf{y}, \mathbf{x}\neq
\mathbf{y}$.

Next, since $f$ is an eigenfunction of $K$, it follows that
\[
\E\{ f(\Yv_k) - f(\Xv_k) \vert\Xv_{k-1}, \Yv_{k-1}
\} =
\lambda\{f(\Yv_{k-1}) - f(\Xv_{k-1})\},
\]
for every $k\geq1$. Therefore,
\begin{eqnarray*}
\Vert K^{n}_{\mathbf{x}^*}-K^{n}_{\mathbf{y}^*}\Vert_{\mathrm{TV}}
&\leq& \E\biggl[ \E\biggl\{\frac{f(\Yv_n)-f( \Xv_n)}{c_1}\Big \vert
\Xv
_{n-1}, \Yv_{n-1}\biggr\}
\Big\vert\Xv_0=\mathbf{x}^*, \Yv_0=\mathbf{y}^* \biggr]\\
&=& \frac{\lambda}{c_1}\E\{f(\Yv_{n-1})-f( \Xv_{n-1})\vert
\Xv
_0=\mathbf{x}^*,
\Yv_0=\mathbf{y}^*\}\\
&=&\frac{\lambda^{n}}{c_1}\{f(\mathbf{y}^*)-f( \mathbf{x}^*)\}.
\end{eqnarray*}
Note that the argument above holds for any $\mathbf{x}^*\preceq
\mathbf{y}^*$.

Note that for any $\mathbf{x}\neq\mathbf{y}$, by the pair-wise dominance
assumption,
there exists
$\zv(\mathbf{x}, \mathbf{y})$ (depends possibly on $\mathbf{x}$ and
$\mathbf{y}$) such that $\zv$
dominates both $\mathbf{x}$
and $\mathbf{y}$ or is dominated by both $\mathbf{x}$ and $\mathbf
{y}$. Hence,
\begin{eqnarray*}
\Vert K^{n}_{\mathbf{x}}- K^{n}_{\mathbf{y}} \Vert_{\mathrm{TV}}
&\leq& \Vert K^{n}_{\mathbf{x}}-K^{n}_{\zv} \Vert_{\mathrm{TV}} +
\Vert
K^{n}_{\mathbf{y}}-K^{n}_{\zv} \Vert_{\mathrm{TV}}\\
&\leq& \frac{\lambda^{n}}{c_1} \vert f(\mathbf{x})- f( \zv) \vert+
\frac
{\lambda^{n}}{c_1} \vert f(\mathbf{y})-f( \zv) \vert\\
&=& \frac{\lambda^{n}}{c_1} \vert f(\mathbf{x}) + f(\mathbf{y}) -
2f( \zv
) \vert.
\end{eqnarray*}
The previous equality follows from the fact that $\zv$ either dominates
or is dominated by
both $\mathbf{x}$ and $\mathbf{y}$, and $f$ is monotone with respect to
$\preceq$,
which implies that
$f(\mathbf{x}) - f(\zv)$ and $f(\mathbf{y}) - f(\zv)$ are either both
positive or
both negative.
Convexity now yields
\[
\Vert K^{n}_{\mathbf{x}}-\pi\Vert_{\mathrm{TV}} \leq\sum_{\mathbf
{y}\in
\mathcal{X}}\pi
(\mathbf{y})\Vert
K^{n}_{\mathbf{x}}- K^{n}_{\mathbf{y}}\Vert_{\mathrm{TV}} \leq\frac
{\lambda^{n}}{c_1}\E
_{\pi} \vert f(\mathbf{x})+
f(\Yv)-2 f(\zv(\mathbf{x},\Yv))\vert.
\]
To get the lower bound, note that
\[
\Vert K^{n}_{\mathbf{x}}-\pi\Vert_{\mathrm{TV}} \geq\frac{1}{2c_2}
\vert\E
_{K^{n}_{\mathbf{x}}}(f(\Yv))
-\E_{\pi}(f(\Yv)) \vert\geq\frac{\lambda^{n}}{2c_2}\vert
f(\mathbf{x})\vert.
\]
Hence the theorem is proved.
\end{pf}

\begin{rem*}
(1) It is to be noted that Theorem~\ref{them1} works for any
arbitrary starting point
without requiring the assumption of reversibility. In Section~\ref
{secmoran}, we show that
the bounds on the total variation distance can be obtained without
explicit knowledge of the
stationary distribution.\vspace*{-6pt}
\begin{longlist}[(1)]
\item[(2)] In all our examples, there will be a unique minimal element
(and no maximal element),
which is clearly sufficient to satisfy the \textit{pair-wise dominance
condition}.
\end{longlist}
\end{rem*}

We now apply this general result for a variety of Markov chains in
population genetics.

\section{Applications}\label{sec3}

\subsection{The Moran process in population genetics} \label{secmoran}

The classical Moran process in population genetics models the evolution
of a population
of constant size by random replacement followed by mutation. Suppose
there are $d$
species in a population of size $N$. At each step, one individual is
chosen uniformly
to die and independently another is chosen uniformly to reproduce.
\textit{They may be
the same individual}. If the latter is of species $i$, the offspring
has probability
$m_{ij}, 1\leq j \leq d$, to mutate to type $j$. Let $\Xv_n = (X_{n1},
\ldots,X_{nd})$
be the vector of counts of species\vadjust{\goodbreak} $1,2, \ldots, d$ at the $n${th}
step. Let
$\mathbb{N}_0 := \mathbb{N} \cup\{0\}$. Then $\{\Xv_n\}_{n \geq0}$
forms a Markov
chain on $\mathcal{X}^d_N$, where
\[
\mathcal{X}_N^d = \Biggl\{\mathbf{x} = (x_1,\ldots,x_d) \in\mathbb{N}_0^d\dvtx
\sum_{i=1}^d x_i = N\Biggr\}.
\]

Let $K$ denote the transition density of this Markov chain. Note that
the size of the state space is $\vert
\mathcal{X}^d_N \vert= {N+d-1 \choose N}$. The one-step transition
probabilities are
%
\begin{eqnarray}\label{mo1}
K(\mathbf{x}, \mathbf{x}+ \ev_i - \ev_j) &=& \frac{x_j}{N}
\Biggl(\sum
_{k=1}^{d} \frac
{x_k}{N} m_{ki}\Biggr), \qquad 1 \leq i \ne j \leq d;
\nonumber\\
K(\mathbf{x}, \mathbf{x}) &=& 1 - \sum_{i\neq j} K(\mathbf{x},
\mathbf{x}+ \ev
_i - \ev_j); \\
K(\mathbf{x},\mathbf{y}) &=& 0 \qquad  \mbox{otherwise},\nonumber
\end{eqnarray}
where $\ev_i$ is the unit vector with $i${th} entry equal to $1$.
\textit{The mutation matrix~$\uM$ is
assumed to be irreducible}. This ensures the irreducibility and
aperiodicity of the transition function
$K$; see proof in the \hyperref[app]{Appendix}. Hence, the stationary distribution of
$K$ exists. Let $\pi$ denote the
density of the stationary distribution with respect to the counting measure.

This model ($d=2$) is due to Moran~\cite{Moran58}. Background and
references can be found in the text by
Ewens~\cite{Ewens04}. When $d=2$, in the continuous-time setting,
Donnelly and Rodrigues~\cite{Donnelly00}
obtain an upper bound in terms of the separation and total variation
distances, when all the individuals
belong to the same generation initially. Watkins~\cite{wat10} analyzes
the infinite allele Moran model in
the discrete-time setting. However, unlike the multi-allele case, the
(infinite) vector of species counts
does not form a Markov chain. Instead, the $N$-dimensional vector whose
$i$th entry is the number of
species with $i$ individuals at the current stage, forms a Markov
chain. It is this fundamentally different
Markov chain that is analyzed in Watkins~\cite{wat10} using strong
stationary times.

In the multi-allele case, which we analyze, a standard choice of the
mutation matrix $\uM= \{m_{ij}\}_{1
\leq i, j \leq d}$ is
\begin{equation}\label{eqnsym-mutation}
\uM= (1-m) \uI+ m \uP,
\end{equation}
where $0 < m \leq1$ is the mutation probability of the offspring, and
$\uP$ is a stochastic matrix with each
row $(p_{1},\ldots,p_{d})$, a probability vector with positive entries.
If mutation happens, the offspring
will change to species $i$ with probability $p_i$. It is known from the
literature that for this standard choice
of the mutation matrix~$\uM$, the corresponding Markov chain is
reversible. Khare and Zhou~\cite{Khare09}
analyze this Markov chain and provide nonasymptotic convergence bounds
in terms of the ``chi-square
distance'' for some natural selected starting points. \textit{In this
paper, we generalize this analysis
in two directions. First, instead of considering the choice $\uM=
(1-m) \uI+ m\uP$,\vadjust{\goodbreak} we consider
a general subclass of mutation matrices described in \textup{(\ref
{eqc1})--(\ref{eqc3})} which includes this
choice as a special case. Second, we provide nonasymptotic convergence
bounds from an arbitrary
starting point}. Consider the class of mutation matrices $\uM$
satisfying one of the monotonicity
conditions specified below:
\begin{equation}\label{eqc1}
m_{dj} < \min_{1 \leq k \leq d-1} m_{kj} \qquad \mbox{for every } 1\leq j
\leq d-1
\end{equation}
or,
\begin{equation}\label{eqc2}
m_{dj} \leq\min_{1\leq k \leq d-1} m_{kj} \qquad \mbox{for every }
1\leq j \leq d-1
\end{equation}
and
$
\uM^* = \{m^*_{ij}\}_{1\leq i,j \leq d-1}$
is irreducible, where $m^*_{ij} = m_{ij} - m_{dj}$
or,
\begin{equation}
m_{dj} \leq\min_{1\leq k \leq d-1}m_{kj}\qquad  \mbox{for every } 1
\leq
j \leq d-1 \label{eqc3}
\end{equation}
and $\uM^*$ has an eigenvector which has all strictly
positive entries.

Each of these conditions essentially says that there is a species,
which we call species $d$ without loss
of generality, such that the mutation probability from this species to
any species is smaller than the
mutation probability from every other species to this species.

It is to be noted that for a general $\uM$ satisfying any one of these
three conditions, the Markov kernel
$K$ is nonreversible, and in this case, often the stationary
distribution of $K$ is not known. Note that condition (\ref{eqc3}) is
satisfied by the standard choice of $\uM= (1-m) \uI+ m\uP$, and
hence the
analysis of this standard choice will come out as a special case. An
example where conditions (\ref{eqc1})
and (\ref{eqc2}) are satisfied would be the following: Suppose $m_{d1}
= \delta$ and $m_{dd} = 1 -
\delta$, that is, the offspring born to species $d$ can possibly mutate
only to species $1$ with a small
probability $\delta$. Suppose $m_{1d} > 0$, that is, species $1$ can
also mutate to species $d$ with a
positive probability. If all the mutation probabilities among species
$1,2, \ldots, d-1$ are larger than
$\delta$, that is, $m_{ij} > \delta$ for $1 \leq i,j \leq d-1$, then
conditions (\ref{eqc1}) and~(\ref{eqc2})
are satisfied.

Let us introduce a partial ordering on $\mathcal{X}_N^{d}$. We define
$\mathbf{x}, \mathbf{y}\in\mathcal{X}_N^{d}$ to
be partially ordered, that is, $\mathbf{x}\preceq\mathbf{y}$ if $x_i
\leq
y_i,  i =
1,2, \ldots, d-1$. This automatically
implies $x_d \geq y_d$. To get bounds on the total variation distance,
according to Theorem~\ref{them1},
we need an eigenfunction $f$ which is strictly monotone in $\preceq$,
that is, if $\mathbf{x}, \mathbf{y}\in
\mathcal{X}^d_{N}$ with $\mathbf{x}\preceq\mathbf{y}$, then
$f(\mathbf{x})\leq f(\mathbf{y})$.
%
\begin{prop} \label{prop1}
Let $K$ denote the Moran process specified by (\ref{mo1}), and
suppose the mutation matrix $\uM$
satisfies any one of conditions (\ref{eqc1})--(\ref{eqc3}). Then $K$
has a linear and strictly monotone
eigenfunction $f$.
\end{prop}

\begin{pf}
Note that
\begin{eqnarray*}
\E_{K(\mathbf{x},\cdot)} [\Xv]
&=& \sum_{1 \le i \ne j \le d} (\mathbf{x}+ \ev_i - \ev_j) \frac
{x_j}{N}
\Biggl( \sum^d_{k=1} \frac{x_k}{N} m_{ki}
\Biggr)\\
&&{}+ \mathbf{x}\Biggl( 1 - \sum_{1 \le i \ne j \le d} \frac{x_j}{N}
\sum
^d_{k=1} \frac{x_k}{N} m_{ki} \Biggr)\\
&=& \mathbf{x}+ \sum_{1 \le i \ne j \le d} (\ev_i - \ev_j) \frac
{x_j}{N}
\Biggl( \sum^d_{k=1} \frac{x_k}{N} m_{ki}
\Biggr)\\
&=& \mathbf{x}+ \sum_{1 \le i, j \le d} (\ev_i - \ev_j) \frac{x_j}{N}
\Biggl(
\sum^d_{k=1} \frac{x_k}{N} m_{ki}
\Biggr)\\
&=& \mathbf{x}+ \sum_{1 \le i \le d} \ev_i \Biggl( \sum^d_{k=1}
\frac{x_k}{N}
m_{ki} \Biggr) - \sum_{1 \le i, j \le
d} \ev_j \frac{x_j}{N} \Biggl( \sum^d_{k=1} \frac{x_k}{N} m_{ki}
\Biggr)\\
&=& \biggl\{ \biggl( 1 - \frac{1}{N} \biggr) \uI_d + \frac{1}{N}
\uM^{T}
\biggr\} \mathbf{x}.
\end{eqnarray*}

Let $\tilde{\av} = (\tilde{a}_i)_{1 \leq i \leq d}$ be any eigenvector
corresponding to an eigenvalue
$\tilde{\lambda}$ of $\uM$. Then we have
\[
\E_{K(\mathbf{x},\cdot)}[\tilde{\av}^T \Xv]
= \biggl\{ \biggl( 1 - \frac{1}{N} \biggr) \tilde{\av}^{T} + \frac
{1}{N}(\uM\tilde{\av})^{T} \biggr\} \mathbf{x}
= \biggl\{ \biggl( 1 - \frac{1}{N} \biggr) + \frac{1}{N} \tilde
{\lambda}
\biggr\} \tilde{\av}^{T} \mathbf{x}.
\]
Hence, $f(\mathbf{x}) = \sum_{i=1}^d \tilde{a}_i x_i$ is an
eigenfunction of
$K$ corresponding to the eigenvalue
$( 1 - \frac{1}{N} ) + \frac{\tilde{\lambda}}{N}$.

We now show that $\uM$ has an eigenvector $\av$ such that $a_i > a_d$
for every $1 \leq i \leq d-1$. It follows
from condition (\ref{eqc1}) that $m^*_{ij}>0$, and from condition
(\ref{eqc2}) that $m^*_{ij} \geq0$, and
$\uM^*$ is irreducible. Hence, under condition (\ref{eqc1}) or (\ref
{eqc2}), by the Perron--Frobenius theorem,
the largest eigenvalue $\lambda^*$ of $\uM^*$ is positive with
multiplicity $1$, and there exists an eigenvector
$\av^* = (a_j^*)_{1 \leq j \leq d-1}$ corresponding to $\lambda^*$,
such that $\av^*$ has all positive entries.
Also, in condition (\ref{eqc3}), we have directly assumed $\av^*$ has
all positive entries. Note that
\begin{eqnarray*}
\lambda^*
&\leq& \max_{1 \leq i \leq d-1} \sum_{j=1}^{d-1} m^*_{ij}\\
&=& \max_{1 \leq i \leq d-1} \sum_{j=1}^{d-1} (m_{ij} - m_{dj})\\
&=& \max_{1 \leq i \leq d-1} (m_{dd} - m_{id})\\
&\leq& m_{dd}\\
&<& 1,
\end{eqnarray*}
since the mutation matrix $\uM$ is assumed to be irreducible.

Let $c$ be defined by
\[
c = \frac{\sum_{j=1}^{d-1} m_{dj} a_j^*}{\lambda^* - 1},
\]
and $\av$ be defined by
\[
a_i= \cases{
a^*_i+c,& \quad$\mbox{if $1\leq i \leq d-1$},$ \vspace*{2pt}\cr
c,& \quad$\mbox{if $i=d$.}$}
\]

Note that, by the definition of $c$,
%
\begin{equation}\label{eqp1}
\sum_{j=1}^d m_{dj} a_j = \sum^{d-1}_{j=1} m_{dj} a_j^* + c =
(\lambda
^*-1)c + c = \lambda^* c.
\end{equation}
We have
\[
\uM^{\ast}\av^{\ast} = \lambda^{\ast} \av^{\ast}\quad\implies\quad\sum^{d-1}_{j=1}
(m_{ij}-m_{dj})(a_j-c) = \lambda^{\ast} (a_i - c)\qquad
\forall  1\leq i\leq d-1.
\]
Note that $\sum_{j=1}^{d-1} (m_{ij} - m_{dj}) = m_{dd} - m_{id}$ and
$a_d = c$. It follows that
%
\begin{equation}\label{eqp2}
\sum_{j=1}^d (m_{ij} - m_{dj}) a_j = \lambda^* (a_i - c).
\end{equation}
Adding (\ref{eqp1}) and (\ref{eqp2}), we get $ \uM\av= \lambda^*
\av$. This shows $\av$ is
an eigenvector of $\uM$ corresponding to eigenvalue $\lambda^*$.\vspace*{1pt}

Thus, $f(\mathbf{x})= \sum_{i=1}^d a_i x_i = \sum_{i=1}^{d-1} (a_i -
a_d) x_i
+ N a_d$, which is
strictly monotone with respect to $\preceq$, is an eigenfunction of $K$
corresponding to the
eigenvalue $\lambda= (1 - \frac{1}{N}) + \frac{\lambda^*}{N}$. Since
$\lambda< 1$, it follows
that $\E_\pi[f(\Xv)] = 0$.
\end{pf}

We now show that for the Moran process, $K$ is monotone with respect to
the partial
ordering, $\preceq$.
%
\begin{thmm} \label{them2}
Let $K$ denote the Moran process specified by (\ref{mo1}), where the
mutation matrix $\uM$
satisfies one of the conditions specified in (\ref{eqc1})--(\ref
{eqc3}). Then $K$ is monotone
with respect to the partial ordering, $\preceq$.
\end{thmm}

\begin{pf}
Consider any $\mathbf{x}\in\mathcal{X}_N^d$ and $\mathbf{y}\in
\mathcal{X}_N^{d}$
with $\mathbf{x}\preceq
\mathbf{y}$. We construct two random vectors $\Xv$ and $\Yv$ such that
$\Xv
\preceq\Yv$ with $\Xv
\sim K(\mathbf{x}, \cdot)$ and $\Yv\sim K(\mathbf{y}, \cdot)$.
This will
immediately
imply that $Kf(\mathbf{x})
\leq Kf(\mathbf{y})$ for any monotone function $f$ and any $\mathbf{x},
\mathbf{y}$ with
$\mathbf{x}\preceq\mathbf{y}$.

Let $\mathbf{x}= (x_1,x_2,\ldots,x_d)$ and $\mathbf{y}=
(y_1,y_2,\ldots
,y_d)$. Then
by assumption $x_i\leq
y_i$ for every $1 \leq i \leq d-1$. We now describe the procedure for
obtaining $\Xv$ and
$\Yv$.

In order to specify the coupling argument, consider two populations
with $N$ individuals
each. Population $1$ has $x_i$ individuals of species $i$, and
population $2$ has $y_i$
individuals of species $i$, for every $1 \leq i \leq d$. We label the
individuals in the
two populations as follows. The individuals of the $i$th species of
population $1$
are labeled from $(\sum_{j=1}^{i}x_{j-1}+1)$ to $\sum_{j=1}^{i}x_j,
i=1, 2, \ldots, d$,
taking $x_0=0$. The labeling of the individuals of population $2$ is
done in the following
way:
\begin{itemize}
\item Note that $x_i \leq y_i$ for $i = 1,2,\ldots,d-1$. For the $i$th
species of
population $2$, where $i=1,2,\ldots,d-1$, we give $x_i$ of the
individuals the exact same
labels as those in species $i$ of population $1$. This leaves $y_i-x_i$
``extra individuals'' to be labeled later.
\item Note that $x_d \geq y_d$. For the $d$th species of population
$2$, the $y_d$
individuals of the $d$th species get exactly same labels as the first
$y_d$ individuals
of the $d$th species of population $1$.
\item Finally, all the $x_d - y_d$ ``extra individuals'' left over in
the first $d-1$
species of population $2$ get the $x_d - y_d$ labels in the $d$th species
of population $1$ which were not assigned in the previous step.
\end{itemize}

The following example illustrates the labeling technique of the $N$
individuals in
population $1$ and population $2$. Consider $N=17$ individuals who
belong to $d=4$
different species type. Also consider $\mathbf{x}=\{1,5,7,4\}$ and
$\mathbf{y}=\{
2,5,8,2\}$. The table
below illustrates the labeling technique.
%
\begin{table}
\caption{Labeling of individuals of population $1$ and population $2$}
\label{tablelabel}
\begin{tabular*}{\textwidth}{@{\extracolsep{\fill}}lcccc@{}}
\hline
\textbf{Species} &\textbf{1st}&\textbf{2nd} &\textbf{3rd} &\textbf{4th}\\
\hline
$\mathbf{x}$&$/$& $/ / / / /$& $/ / / / / / /$&$/ / / /$\\
labels&$1$&$2,3,4,5, 6$&$7,8,9,10,11,12,13$&$14,15,{\bf16},{\bf17}$\\
$\mathbf{y}$&$/ /$& $/ / / / /$&$/ / / / / / / /$&$/ /$\\
labels&$1,{\bf16}$&$2,3,4,5,6$&$7, 8,9,10,11,12,13,{\bf17}$&$14,15$\\
\hline
\end{tabular*}
\end{table}

In Table~\ref{tablelabel}, we label the individuals of population $1$
from $1$ to $17$ based
on $\mathbf{x}$. For the $1$st species of population $2$, there are $2$
individuals, the first
individual gets the label $1$, same as the label of the first
individual of population $1$,
and the second individual is an ``extra individual,'' to be labeled
later. Now, for the 2{nd} species, there are the same number of individuals for both the
populations, so these
individuals get the same labels. For the $3$rd species, there is one
``extra individual,''
to be labeled later; other individuals get the same labels. The $4$th
species has $2$
individuals in population $2$, who get the same labels as the first $2$
individuals of the
$4$th species in population $1$. Last, $2$ extra labels $16$ and $17$
are assigned to
the ``extra individuals'' of species $1$ and $3$ of population $2$,
respectively.

Let us return to the general proof, and define $k_1 := \sum_{i=1}^{d-1}
x_i$ to be the total
number of individuals in the first $d-1$ species of population $1$ and
$k_2 := y_d$ to be
the number of individuals in species $d$ of population $2$. We now
change the species
configuration of population $1$ and population $2$ in four sub-steps
which are described
below:
\begin{longlist}[(III)]
\item[(I)] Choose a label uniformly between $1$ to $N$. Call it $i_1$.
\item[(II)] Independently choose another label uniformly between $1$ to
$N$. Call it $i_2$.
\item[(III)] Let $s_{1,i_2}$ and $s_{2, i_2}$ denote the species of the
individual
labeled $i_2$ in population $1$ and population $2$, respectively. Add
one individual of
species $s_{1, i_2}$ to population $1$ and one individual of species
$s_{2, i_2}$ to
population $2$.

Note that if $1\leq i_2 \leq k_1 +k_2$, then $s_{1, i_2} = s_{2,
i_2}:=s_{i_2}$. In this
case the newly added individual in both the populations mutates in the
following way:
Generate $U\sim\operatorname{Uniform}[0,1]$. If $0\leq U < m_{s_{i_2}1}$, the added
individual mutates to species $1$. If $m_{s_{i_2}1}\leq U
<m_{s_{i_2}1}+ m_{s_{i_2}2}$, the
added individual mutates to species $2$, and so on. Finally, if
$m_{s_{i_2}1}+m_{s_{i_2}2}+\cdots+m_{s_{i_2}(d-1)}\leq U \leq1$, the
added individual mutates
to species $d$. Hence, after the mutation, both populations have an
individual of the same
species added, which therefore preserves the partial ordering between
their species
configurations.

Next, suppose $k_1+k_2+1\leq i_2 \leq N$, then $s_{1, i_2}=d$ and
$s_{2, i_2}$ is one of
the first $d-1$ species. Note that $m_{s_{2, i_2}j} \geq m_{dj}$ for
every $j=1,2,\ldots,d-1$.
The newly added individual in population $1$ mutates in the following way:
Generate $U \sim\operatorname{Uniform}[0,1]$. If $0\leq U < m_{d1}$, the added
individual mutates to
species~$1$. If $m_{d1}\leq U <m_{d1} + m_{d2}$, the added individual
mutates to species
$2$, and so on. Finally, if $m_{d1}+m_{d2}+\cdots+m_{d(d-1)}\leq U
\leq
1$, the added
individual mutates to species $d$. Now, in population $2$, the newly
added individual mutates
in the following way: Choose the same $U$ as for population $1$.
If $0\leq U < m_{d1}$ or $m_{d1} + m_{d2} + \cdots+
m_{d(d-1)} \leq U <
m_{s_{2, i_2}1} + m_{d2} + \cdots+ m_{d(d-1)}$,
the individual
mutates to species~$1$. If $m_{d1} \leq U < m_{d1} + m_{d2}$ or $m_{s_{2, i_2}1} + m_{d2} +
\cdots+ m_{d(d-1)} \leq U < m_{s_{2, i_2}1} + m_{s_{2, i_2}2} + m_{d3}
+ \cdots+
m_{d(d-1)}$, the individual mutates to species $2$, and so on. Finally,
if $m_{s_{2, i_2}1} + m_{s_{2, i_2}2} + \cdots+ m_{s_{2, i_2}(d-1)} \leq U
\leq1$, the
individual mutates to species $d$. Hence, when $0 \leq U \leq m_{d1}
+ m_{d2} + \cdots+
m_{d(d-1)}$ or when $m_{s_{2, i_2}1} + m_{s_{2, i_2}2} + \cdots+
m_{s_{2, i_2}(d-1)} \leq U
\leq1$, the newly added individual in both the populations mutate to
the same species, which
preserves the partial ordering between their species configurations.
Alternatively, if
$m_{d1} + m_{d2} + \cdots+ m_{d(d-1)} \leq U \leq m_{s_{2, i_2}1} +
m_{s_{2, i_2}2} + \cdots
+ m_{s_{2, i_2}(d-1)}$, then after mutation the newly added individual
in the population $1$
is in species $d$, but the newly added individual in the population $2$
is in any of the
first $d-1$ species. This again preserves the partial ordering between
the species
configurations in population $1$ and population~$2$.

\item[(IV)] Finally, the individual corresponding to the label $i_1$
dies for both the
populations. If $1 \leq i_1 \leq k_1 + k_2$, then the individual
belongs to the same
species for both the populations. If $k_1 + k_2 + 1 \leq i_1 \leq N$,
then the individual
corresponding to the label $i_1$ belongs to species $d$ for population
$1$ and is an
``extra individual'' in the first $d-1$ species of population $2$. In
either case, the
partial ordering is preserved.
\end{longlist}

Let $\Xv$ and $\Yv$ be the resulting species configurations of
population $1$ and
population $2$, respectively. Note that marginally the movement from
both $\mathbf{x}$ to $\Xv$
and $\mathbf{y}$ to $\Yv$ follows the transition mechanism of $K$, and
$\Xv
\preceq\Yv$. This
completes the proof.
\end{pf}

\subsubsection{Bounds on total variation distance}\label{sec3.1.1}

For the partial ordering, $\preceq$, discussed above, applying Theorem
\ref{them1} in the
case of the Moran model, provides us with bounds on the total variation
distance. We
have shown that for the Moran process, $K$ is monotone with respect to
the partial
ordering, $\preceq$; see Theorem~\ref{them2}. It is easily seen that
$\mathbf{0}$ (with
first $d-1$ entries equal to zero, and the $d$th entry equal to $N$) is
dominated by
$\mathbf{x}$ for every $\mathbf{x}\in\mathcal{X}^d_{N}$. Hence, the
pair-wise
dominance property is
satisfied. Recall that by Proposition~\ref{prop1}, there exists an
eigenfunction $f(\mathbf{x})
= \sum^d_{i=1} a_i x_i = \sum_{i=1}^{d-1} (a_i - a_d) x_i + N a_d$ of
$K$ corresponding
to the eigenvalue $\lambda= 1 - \frac{1}{N} + \frac{\lambda^*}{N}$,
such that $f$ is
strictly monotone with respect to the partial ordering, $\preceq$.
Hence, the conditions
of Theorem~\ref{them1} are satisfied, and the bounds on the total
variation distance are
obtained as
\begin{eqnarray*}
& & \frac{\lambda^{n}}{2 c_2} \vert f(\mathbf{x}) \vert\leq\Vert
K^{n}_{\mathbf{x}}
- \pi
\Vert_{\mathrm{TV}} \leq\frac{\lambda^{n}}{c_1} \E_\pi\{f(\Yv) +
f(\mathbf{x})
- 2
f(\mathbf{0})\}\\
&&\qquad\implies\quad \frac{\lambda^{n}}{2 c_2} \vert f(\mathbf{x}) \vert\leq
\Vert
K^{n}_{\mathbf{x}} - \pi
\Vert_{\mathrm{TV}} \leq\frac{\lambda^{n}}{c_1} \{ f(\mathbf{x}) - 2
f(\mathbf{0})\}
\\
&&\qquad\implies\quad \frac{\lambda^{n}}{2 c_2} \Biggl\vert\sum_{i=1}^{d-1} a^*_i
x_i + N a_d
\Biggr\vert\leq\Vert K^{n}_{\mathbf{x}} - \pi\Vert_{\mathrm{TV}}
\leq
\frac
{\lambda^{n}}{c_1}
\Biggl\{ \sum_{i=1}^{d-1} a^*_i x_i - N a_d \Biggr\},
\end{eqnarray*}
where $a_i^* = a_i - a_d > 0$ for every $1 \leq i \leq d - 1$ (by the
monotonicity of
$f$), $c_1 =  \min_{1 \leq i \leq d-1} a^*_i > 0$ and $c_2
= \max\{
-N a_d, N (\max_{1 \leq i \leq d - 1} a_i^* + a_d )
\}
$. Note that $a_d
< 0$. Note again that the stationary distribution $\pi$ in not known in
general, but the
analysis above leads to upper and lower bounds which do not depend on
the stationary
distribution, and are reasonably close to each other.

\subsubsection{Bounds on total variation distance in the special case}
\label{secspmoran}
We now provide a nonasymptotic convergence analysis for the special
choice of $\uM=
(1-m) \uI+ m \uP$. It has been proved earlier in Khare and Zhou \cite
{Khare09}
that the Markov chain $K$ corresponding to the multi-allele Moran model
with $\uM=
(1-m) \uI+ m \uP$ has second largest eigenvalue $\lambda= 1 - \frac
{\vert\ualpha
\vert}{N (N + \vert\ualpha\vert)}$, where $\vert\ualpha\vert:=
\sum
^{d}_{i=1}
\alpha_i$, $\ualpha= (\alpha_1, \alpha_2, \ldots, \alpha_d)$, where
$\alpha_i =
\frac{Nmp_i}{1-m}$, with the eigenspace given by the space of centered
linear functions
of $x_1, x_2,\ldots, x_{d-1}$. After simplification, we obtain
$\lambda
= 1-\frac{m}{N}$.
It is known that the stationary distribution in this case is the
Dirichlet-multinomial
distribution with parameters $N$ and $\ualpha$. The
Dirichlet-multinomial distribution
with parameters $N > 0$ and $\bolds{\alpha} = (\alpha_1, \alpha_2,
\ldots,
\alpha_d), \alpha_i>0$, has probability mass function given
by
\[
\mathcal{DM}(\mathbf{x}|N,\bolds{\alpha}) =\frac{\prod_{i=1}^d
{x_i + \alpha_i
- 1\choose x_i}}{{N + |\bolds{\alpha}| - 1\choose N}} ,\qquad
\mathbf{x} \in
\mathcal{X}_N^{d}.
\]
Since $\uM= (1-m) \uI+ m \uP$, it follows that $\uM^* = (1-m) \uI
_{d-1}$. Hence,
any $(d-1)$-dimensional vector with positive entries is an eigenvector
of $\uM^*$.
Suppose we choose the eigenvector $\mathbf{a}^*$ of $\uM^*$ such that
$a_i^* = 1$ for
$i < d$. Then, for the Markov chain $K$, we get the eigenfunction
$f(\mathbf{x}
) =
\sum_{i=1}^{d-1} x_i - N(1 - p_d)$ corresponding to the eigenvalue
$\lambda= 1 -
\frac{m}{N}$. Note that $f$ is strictly monotone with respect to the
partial ordering,
$\preceq$. As in the case of the general multi-allele Moran model, here
also it is easily
seen that $\mathbf{0}$ is dominated by $\mathbf{x}$ for every
$\mathbf{x}\in
\mathcal
{X}^d_{N}$. We have
$c_1 = 1$ and $c_2 = \max\{Np_d, N(1-p_d)\}$. Thus, bounds on total
variation distance
are obtained as
%
\begin{eqnarray}\label{eqmo2}
\label{eqmo1}
\Vert K^{n}_{\mathbf{x}} - \pi\Vert_{\mathrm{TV}} &\geq& \frac{ (1 -
{m}/{N})^{n}}{2 \max
\{Np_d, N(1 - p_d)\}} \Biggl\vert\sum_{i=1}^{d-1} x_i - N(1 - p_d)
\Biggr\vert,\\
\Vert K^{n}_{\mathbf{x}} - \pi\Vert_{\mathrm{TV}} &\leq& \biggl( 1 -
\frac{m}{N}
\biggr)^{n}
\Biggl( \sum_{i=1}^{d-1} x_i + N(1 - p_d) \Biggr).
\end{eqnarray}

\begin{eg}
Consider the multi-allele Moran model in the special case when the
mutation matrix
$\uM= (1-m) \uI+ m \uP$. Suppose the population size \mbox{$N = 100$}, with
$d = 5$
species and mutation probability $m = 0.7$. When mutation occurs, the individual
mutates to the $i$th species with probability $p_i = 1/5$. Using~(\ref{eqmo1})
and (\ref{eqmo2}), for a starting state $\mathbf{x}=
(0,10,0,10,80)$, the
bounds on the
total variation distance are obtained as
%
\begin{equation}\label{eqmo3}
0.375 \biggl(1-\frac{7}{1000}\biggr)^{n} \leq\| K_\mathbf{x}^n -
\pi\|
_{\mathrm{TV}}
\leq100 \biggl( 1 - \frac{7}{1000} \biggr)^{n}.
\end{equation}
For $\varepsilon=0.01$, (\ref{eqmo3}) tells us that $516$ steps are
necessary and
$1312$ steps are sufficient for the total variation distance to be less
then $0.01$.
The crude upper bound for the total variation distance is $(
2.1665 \times
10^{15} ) ( 1 - \frac{7}{1000} )^{n}$, which gives
$5683$ steps are
sufficient for the total variation distance to be less then $0.01$.
\end{eg}

\subsection{Sequential P{\'o}lya urn models}\label{polya1}

Choose $d$ urns with $N$ balls distributed in them. Suppose the
inherent weight of urn $i$ is
$\alpha_i$, $i=1,2,\ldots,d$, and let $\ualpha= (\alpha_1, \alpha_2,
\ldots, \alpha_d)$ denote
the vector of urn weights and $\vert\ualpha\vert= \sum_{i=1}^d
\alpha
_i$ denote the total
inherent weight of $d$ urns. Suppose that each ball has unit weight.\looseness=-1\vadjust{\goodbreak}
\begin{longlist}[(1)]
\item[(1)] P\'olya level model~\cite{Khare09}: Consider the Markov
chain whose one-step movement
consists of the following sub-steps:
\begin{enumerate}[(iii)]
\item[(i)] Randomly choose $s$ balls out of $N$ balls and mark them.
\item[(ii)] Draw an urn with probability proportional to its weight
(inherent\break weight${} + {}$weight
of balls) and add a ball (of unit weight) to the chosen urn. Repeat
this~$s$ times.
\item[(iii)] Remove the $s$ marked balls from the respective urns.
\end{enumerate}
\item[(2)] P\'olya up--down model~\cite{Khare09}: These are variations
of P\'olya level models,
where the three steps are performed in the following order (ii), (i)
(with $N+s$ total balls)
and (iii).
\item[(3)] P\'olya down--up model~\cite{Khare09}: These are variations
of P\'olya level models,
where the three steps are performed in the following order (i), (iii)
and (ii).
\end{longlist}

We first analyze the Markov chain corresponding to the P\'olya level
model. Let $X_{ni}$
denote the number of balls in the $i$th urn at the $n$th step of the
P\'olya level
model. Then $\{\Xv_n = (X_{n1}, X_{n2}, \ldots,X_{nd}), n=0,1,2,
\ldots
\}$ forms a multivariate
Markov chain on $\mathcal{X}_N^d$. Let $K$ denote the transition
density of this Markov chain.
Let $\preceq$ be the partial ordering on $\mathcal{X}_N^d$ as in the
multi-allele Moran model.
%
\begin{thmm} \label{them3}
$K$ is monotone with respect to the partial ordering, $\preceq$.
\end{thmm}

\begin{pf}
Consider any $\mathbf{x}\in\mathcal{X}_N^d$ and $\mathbf{y}\in
\mathcal{X}_N^{d}$
with $\mathbf{x}\preceq
\mathbf{y}$. We construct two random vectors $\Xv$ and $\Yv$ such that
$\Xv
\preceq\Yv$ with $\Xv
\sim K(\mathbf{x}, \cdot)$ and $\Yv\sim K(\mathbf{y}, \cdot)$.
This will
immediately
imply that $Kf(\mathbf{x})
\leq Kf(\mathbf{y})$ for any monotone function $f$ and any $\mathbf{x},
\mathbf{y}$ with
$\mathbf{x}\preceq\mathbf{y}$.

In order to specify the coupling argument, we consider two populations
of $N$ balls each,
with $N$ balls distributed in $d$ urns based on $\mathbf{x}$ and
$\mathbf{y}$,
respectively. We use the
same labeling technique for both the populations as discussed in
Theorem~\ref{them2}
(regarding species as urns and individuals as balls).



We now change the urn configuration of population $1$ and population
$2$ in three sub-steps
which are described below:
\begin{longlist}[(III)]
\item[(I)] Choose $s$ labels without replacement from $1$ to $N$.
\item[(II)] This sub-step will consist of $s$ sequential urn draws, and
after each draw, an
extra ball will be added to the chosen urn for both the populations as
described below.
Repeat the following for $j=1, 2, \ldots, s$.

Generate $U_j \sim\operatorname{Uniform}[0,1]$. Now, at the beginning of the
$j$th draw in
this sub-step, there are, in total, $N+j-1$ balls each in both the
populations. Hence the
total weight of the urns (with balls) in both the populations is $\vert
\ualpha\vert+ N +
j - 1$. Let $\Xv^{j-1} := (x^{j-1}_1, x^{j-1}_2, \ldots, x^{j-1}_d)$ be
the configuration of
the balls in the $d$ urns of population $1$ at the beginning of the
$j$th draw, and
$\Yv^{j-1} := (y^{j-1}_1, y^{j-1}_2, \ldots, y^{j-1}_d)$ be the
configuration of
the balls in the $d$ urns of population $2$ at the beginning of the
$j$th draw. Let us
denote the normalized probability vector of the urn weights for
population $1$ by
$\pv^{j-1} = (p^{j-1}_1, p^{j-1}_2, \ldots, p^{j-1}_d)$, and the
normalized probability
vector of urn weights for population $2$ by $\qv^{j-1} = (q^{j-1}_1,
q^{j-1}_2, \ldots,
q^{j-1}_d)$, where $p^{j-1}_i = \frac{\alpha_i + x^{j-1}_i}{\vert
\ualpha\vert+ N + j -
1}$ and $q^{j-1}_i = \frac{\alpha_i + y^{j-1}_i}{\vert\ualpha\vert+
N + j - 1}$.

Procedure to choose an urn for population $1$ at the $j$th draw:

If $0\leq U_j <
p^{j-1}_1$, choose urn $1$. If $p^{j-1}_1 \leq U_j < p^{j-1}_1 +
p^{j-1}_2$, choose urn $2$
and so on. Finally, if $p^{j-1}_1 + p^{j-1}_2 + \cdots+
p^{j-1}_{(d-1)} \leq U_j \leq1$,
choose urn $d$. Add a ball to the chosen urn.

The following is the procedure to choose urn for population $2$ at the
$j$th draw:

$\mbox{If } 0 \leq U_j < p^{j-1}_1 \mbox{ or, }
p^{j-1}_1 + p^{j-1}_2 + \cdots+ p^{j-1}_{(d-1)} \leq U_j <
q^{j-1}_1 + p^{j-1}_2 + \cdots+ p^{j-1}_{(d-1)}, \mbox{ choose urn $1$. }$
$\mbox{If } p^{j-1}_1 \leq U_j <p^{j-1}_1 + p^{j-1}_2  \mbox{ or, }
q^{j-1}_1 + p^{j-1}_2 + \cdots+\break p^{j-1}_{(d-1)} \leq U_j <
q^{j-1}_1 + q^{j-1}_2 + \cdots+ p^{j-1}_{(d-1)}, \mbox{ choose urn $2$}$
and so on. Finally, $\mbox{if } q^{j-1}_1 + q^{j-1}_2 + \cdots+
q^{j-1}_{(d-1)} \leq
U_j \leq1$, choose urn $d$. Add a ball to the chosen urn.
\item[(III)] Remove the balls corresponding to the $s$ labels in
sub-step (I) from both
the populations.
\end{longlist}

It is to be noted that in sub-step (II), assuming $\Xv^{j-1} \preceq
\Yv
^{j-1}$ (and hence
$\pv^{j-1} \preceq\qv^{j-1}$), the mechanism for drawing urns is such
that either the
same urn is chosen for both the populations or when the $d$th urn is
chosen for
population $1$, then any of the first $d-1$ urns is chosen for
population $2$. Hence,
$\Xv^j \preceq\Yv^j$ (and hence $\pv^j \preceq\qv^j$). Since $\Xv
^0 =
\mathbf{x}$ and $\Yv^0 =
\mathbf{y}$, it follows by induction (on $j$) that $\Xv^j \preceq\Yv
^j$ for
$j = 1,2, \ldots,
s$. In sub-step (III), the balls with the same $s$ labels are removed
from both the
populations. Based on the labeling procedure, either balls with the
same label lie in the
same urn for both the populations, or the ball lies in the $d$th urn
for population $1$
and is an ``extra ball'' in the first $d-1$ urns for population $2$. In
either case,
removing balls with the same label from both the populations does not
change the partial
ordering of the urn configurations.

Let $\Xv$ and $\Yv$ be the resulting urn configurations of population
$1$ and population $2$,
respectively. It follows from the discussion above that $\Xv\preceq
\Yv
$. Note that
marginally the movement from both $\mathbf{x}$ to $\Xv$ and $\mathbf{y}$
to $\Yv$
follows the transition
mechanism of $K$. To see this, note that the probability of choosing
the $i$th urn at the
$j$th draw in sub-step (II) for population $1$ is
$\Po(\sum^{i-1}_{\ell=1} p^j_\ell\leq U_j \leq\sum^{i}_{\ell=1}
p^j_\ell) = p^j_i$; and
the corresponding probability for population $2$ is
$\Po(\sum^{i-1}_{\ell=1} p^j_\ell\leq U_j \leq\sum^{i}_{\ell=1}
p^j_\ell) + \Po
(\sum^{i-1}_{\ell=1} q^j_\ell+ \sum^{d-1}_{\ell=i} p^j_\ell\leq U_j
\leq
\sum^{i}_{\ell=1} q^j_\ell+ \sum^{d-1}_{\ell=i+1} p^j_\ell) = q^j_i$.
This completes the
proof.
\end{pf}

We can similarly argue that the Markov chain corresponding to the P{\'
o}lya up--down model and
the P{\'o}lya down--up model are stochastically monotone with
respect to the partial ordering, $\preceq$ in $\mathcal{X}^d_N$.

\subsubsection{Bounds on total variation distance} \label{polyabound}

In case of the P\'{o}lya level model, the second largest eigenvalue
$\lambda= 1 - \frac{s \vert
\ualpha\vert}{N (N + \vert\ualpha\vert)}$. We know that the stationary
distribution of the
P\'{o}lya level model is the Dirichlet-multinomial distribution with
parameters $N$ and
$\ualpha$. The eigenfunction $f(\mathbf{x}) = \sum_{i=1}^{d-1} x_i - N
(
1 - \frac{\alpha_d}{\vert
\ualpha\vert} )$ corresponding to $\lambda$ is strictly monotone
in $\preceq$. Let
$\mathbf{0}$ be a $d$-dimensional vector such that the first $d-1$ entires
are zero, and the $d$th
entry is $N$. It is easily seen that $\mathbf{0}$ is dominated by
$\mathbf{x}$ for
every $\mathbf{x}\in
\mathcal{X}^{d}_{N}$. Hence, the conditions of Theorem~\ref{them1} are
satisfied, with $c_1 = 1$
and $c_2 = \max\{N\frac{\alpha_d}{\vert\ualpha\vert},
N(1
-\frac{\alpha_d}{\vert\ualpha\vert})\}$. Let
$p_i=\frac
{\alpha_i}{\vert
\ualpha\vert}, i=1,2,\ldots,d$.
Thus, the bounds on the total variation distance are obtained as
%
\begin{eqnarray}
\label{eqpol1}
\Vert K^{n}_{\mathbf{x}} - \pi\Vert_{\mathrm{TV}} &\geq& \frac
{\lambda^{n}}{2
\max\{Np_d, N(1 - p_d)\}} \Biggl\vert\sum_{i=1}^{d-1}
x_i -
N(1 - p_d) \Biggr\vert,\\
\label{eqpol2}
\Vert K^{n}_{\mathbf{x}} - \pi\Vert_{\mathrm{TV}} &\leq& \lambda^{n}
\Biggl( \sum
_{i=1}^{d-1} x_i + N
(1 -p_d) \Biggr).
\end{eqnarray}
Similarly, in the case of P{\'o}lya down--up models, the second largest
eigenvalue is given by
$\lambda= (1 - \frac{s}{N}) (1 - \frac{s}{N + \vert\ualpha
\vert})^{-1}$ and in the
case of P{\'o}lya up--down models, the second largest eigenvalue is
given by $\lambda=
(1 + \frac{s}{N})^{-1}(1 + \frac{s}{N + \vert\ualpha\vert})$. These
can be substituted in
(\ref{eqpol1}) and (\ref{eqpol2}) to get the corresponding total
variation bounds for these
models.
%
\begin{rem*}
Note that the coefficient of $\lambda^n$ in the upper bound derived in
(\ref{eqpol2}) is
at most $2N$. Let us try and compare it to the coefficient of $\lambda
^n$ in the crude upper
bound, which is given by
\[
\frac{1}{2 \sqrt{\pi(\mathbf{x})}} = \frac{1}{2}\sqrt{\frac{{{N +
|\bolds{\ualpha}| - 1 \choose
N}}}{\prod_{i=1}^d {{x_i + \alpha_i -1 \choose x_i}}}}.
\]
At one possible extreme, when all entries of $\mathbf{x}$ except the
$i$th one
are zero,\vspace*{1pt} the
coefficient is essentially a polynomial in $N$ of degree $\frac{\vert
\ualpha\vert-
\alpha_i}{2}$. At the other possible extreme, when all the entries of
$\mathbf{x}$ are equal to
$\frac{N}{d}$ (assuming $\frac{N}{d}$ is an integer), the coefficient
is essentially a
polynomial in $N$ of degree $\frac{d - 1}{2}$. The main fact is that
the~coefficient of
$\lambda^n$ in the upper bound derived in (\ref{eqpol2}) is linear in
$N$, whereas the
coefficient of $\lambda^n$ in the crude upper bound almost always
behaves like a polynomial
of a higher degree in $N$.
\end{rem*}

\begin{eg}
Consider the P\'{o}lya level model where $N=100$ balls are distributed
in $d=5$ urns.
Suppose $s=2$ balls are chosen and each urn has inherent weight $\alpha
_i = 180$ for
every $1 \leq i \leq5$. Using (\ref{eqpol1}) and (\ref{eqpol2}), for
a starting
state $\mathbf{x}=(0,20,0,20,60)$, the bounds on the total variation distance
are obtained
as
%
\begin{equation}\label{eqpol3}
0.25 \biggl( 1 - \frac{9}{500} \biggr)^{n} \leq\| K_\mathbf{x}^n -
\pi
\|_{\mathrm{TV}}
\leq120 \biggl( 1 - \frac{9}{500} \biggr)^{n}.
\end{equation}
For $\varepsilon=0.01$, (\ref{eqpol3}) tells us that $178$ steps are
necessary and
$518$ steps are sufficient for the total variation distance to be less
than $0.01$.
The crude upper bound for total variation distance is $( 6.1094
\times10^{13}
) ( 1 - \frac{9}{500} )^n$ which would have implied
$2002$ steps
are sufficient for the total variation distance to be less then $0.01$.
\end{eg}

\subsection{A generalized Ehrenfest urn model} \label{secEhrenfest}

There are $N$ indistinguishable balls to be distributed to $d$ urns. At
each step, $s$ balls
are chosen at random from the total of $N$ balls, and each of them is
redistributed
independently according to the same probability $\pv=(p_1, p_2,\ldots
,p_d)$. Let $X_{ni}$ be the
number of balls in the $i$th urn at the $n$th step of the Markov chain.
Then $\{\Xv_{n}
= (X_{n1}, X_{n2}, \ldots, X_{nd}),  n = 0,1,2, \ldots\}$ forms a
multivariate Markov chain
on $\mathcal{X}_N^{d}$. Let $K$ denote the transition density of this
Markov chain.

Consider the same partial ordering, $\preceq$, as defined in the case
of the Moran process.
We now show that $K$ is a monotone Markov chain with respect to the
partial ordering,
$\preceq$.
%
\begin{thmm} \label{them4}
$K$ is monotone with respect to the partial ordering, $\preceq$.
\end{thmm}

\begin{pf}
Consider any $\mathbf{x}\in\mathcal{X}_N^d$ and $\mathbf{y}\in
\mathcal{X}_N^{d}$
with $\mathbf{x}\preceq
\mathbf{y}$. We construct two random vectors $\Xv$ and $\Yv$ such that
$\Xv
\preceq\Yv$ with $\Xv
\sim K(\mathbf{x}, \cdot)$ and $\Yv\sim K(\mathbf{y}, \cdot)$.
This will
immediately
imply that $Kf(\mathbf{x})
\leq Kf(\mathbf{y})$ for any monotone function $f$ and any $\mathbf{x},
\mathbf{y}$ with
$\mathbf{x}\preceq\mathbf{y}$.

In order to specify the coupling argument, we consider two populations
of $N$ balls each,
with $N$ balls distributed in $d$ urns based on $\mathbf{x}$ and
$\mathbf{y}$,
respectively. We use the
same labeling technique for both the populations as discussed in
Theorem~\ref{them2}
(regarding species as urns and individuals as balls).

We now change the urn configuration of population $1$ and population
$2$ in five
sub-steps which are described below:
\begin{longlist}[(III)]
\item[(I)] Choose $s$ labels without replacement from $1$ to $N$.
\item[(II)] Remove the balls with the chosen labels from both $\mathbf
{x}$ and
$\mathbf{y}$.
\item[(III)] Choose an urn, such that urn $i$ is chosen with
probability $p_i$ for every $i =
1,2, \ldots, d$.
\item[(IV)] Add a ball to the chosen urn for both the current $\Xv$ and
$\Yv$ configurations.
\item[(V)] Repeat steps (III) and (IV) $s$ times independently.
\end{longlist}

Let $k_1 := \sum_{i=1}^{d-1} x_i$ be the total number of balls in the
first $d-1$ urns of
population $1$, and $k_2 := y_d$ be the number of balls in the $d$th
urn of
population~$2$. Consider sub-steps (I) and (II). Without loss of
generality, let us assume
out of $s$ labels chosen, $r$ labels are between $1$ and $k_1+k_2$ and
$s-r$ labels are
between $k_1+k_2+1$ and $N$.
\begin{itemize}
\item Each of the $r$ balls corresponding to labels $1$ to $k_1+k_2$
lies in exactly the
same urn for both the populations. Removing these does not change the
partial ordering between
the urn configurations.
\item Each of the $s-r$ balls corresponding to labels $k_1+k_2+1$ to
$N$ lie in urn $d$ for
population $1$, and are ``extra balls'' lying in the first $d - 1$ urns
for population~$2$.
Hence, removing them does not change the partial ordering between the
urn configurations of
population $1$ and population $2$.
\end{itemize}

Consider sub-steps (III), (IV) and (V). Since the balls are put in the
same urn for both
the populations, adding the new balls does not change the partial
ordering between the urn
configurations.

Let $\Xv$ and $\Yv$ be the resulting urn configurations of population
$1$ and
population~$2$, respectively. Note that marginally the movement from
both $\mathbf{x}$ to $\Xv$ and
$\mathbf{y}$ to $\Yv$ follows the transition mechanism of $K$, and
$\Xv
\preceq\Yv$. This
completes the proof.
\end{pf}

The following example illustrates the one-step movement of the above
construction in
population $1$ and population $2$ for Theorem~\ref{them4}.
\begin{eg}
Consider the same $\mathbf{x}$ and $\mathbf{y}$ as in Table \ref
{tablelabel}.
Suppose the $4$ balls
chosen in sub-step (I) are with labels $6$, $8$, $14$ and $16$. It is
evident that the
removal of the balls with the chosen labels in sub-step (II) does not
alter the partial
ordering between the urn configurations of the two populations. Since
the urn chosen
in sub-step (III) is same for both the populations, adding a ball to
the urn in
sub-step (IV) does not change the partial ordering between the urn
configurations of
the two populations.
\end{eg}

\subsubsection{Bounds on total variation distance} \label{ehrenbound}

For the partial ordering, $\preceq$, discussed above, applying Theorem
\ref{them1} in the case
of the generalized Ehrenfest urn model, provides us with bounds on the
total variation
distance.

It has been proved earlier in Khare and Zhou~\cite{Khare09} that the
generalized Ehrenfest
urn model has second largest eigenvalue $\lambda= 1 - \frac{s}{N}$,
with the eigenspace given
by the space of linear functions of $x_1, x_2, \ldots, x_{d-1}$. It is
known that the
stationary distribution is the multinomial distribution with parameters
$N$ and $\pv$. The
eigenfunction $f(\mathbf{x}) = p_d \sum_{i=1}^{d-1} x_i - (1 - p_d) x_d
= \sum
_{i=1}^{d-1} x_i - N(1
- p_d)$ corresponding to the eigenvalue $\lambda$ is strictly monotone
in $\preceq$. Again, it
is easily seen that $\mathbf{0}$ is dominated by $\mathbf{x}$, for every
$\mathbf{x}\in
\mathcal{X}^d_N$. Hence,
the conditions of Theorem~\ref{them1} are satisfied. We have $c_1 = 1$
and $c_2 = \max\{Np_d,
N(1-p_d)\}$. Thus, the bounds on total variation distance are
%
\begin{eqnarray}
\label{eqehf1}
\Vert K^{n}_{\mathbf{x}} - \pi\Vert_{\mathrm{TV}} &\geq& \frac{(1 -
{s}/{N})^{n}}{2 \max\{Np_d,
N(1 - p_d)\}} \Biggl\vert\sum_{i=1}^{d-1} x_i - N(1 - p_d)
\Biggr\vert,\\
\label{eqehf2}
\Vert K^{n}_{\mathbf{x}} - \pi\Vert_{\mathrm{TV}} &\leq& \biggl( 1 -
\frac{s}{N}
\biggr)^{n} \biggl(
\sum_{i=1}^{d-1} x_i + N(1 - p_d) \biggr).
\end{eqnarray}

\begin{rem*}
Note that the coefficient of $( 1 - \frac{s}{N} )^n$ in the
upper bound derived in
(\ref{eqehf2}) is at most $2N$. We compare it to the coefficient of
$( 1 - \frac{s}{N}
)^n$ in the crude upper bound, which is given by
\[
\frac{1}{2 \sqrt{\pi(\mathbf{x})}} = \frac{1}{2 \sqrt{{N\choose
{\mathbf{x}}}}} \prod
_{i=1}^d \biggl(
\frac{1}{\sqrt{p_i}} \biggr)^{x_i}.
\]
At one possible extreme, when all entries of $\mathbf{x}$ except the
$i$th one
are zero, the
coefficient is $\frac{1}{2} ( \frac{1}{\sqrt{p_i}}
)^N$. At
the other possible
extreme, when all the entries of $\mathbf{x}$ are equal to $\frac{N}{d}$
(assuming $\frac{N}{d}$ is
an integer), using Stirling's approximation for large~$N$,\footnote{Note that $N$ is the
notation for the total number of balls in the urns, not the number of
steps.} the
coefficient is
\[
\frac{(2 \pi N)^{{(d-1)}/{4}}}{2 d^{{d}/{4}}} \biggl( \frac
{1}{\sqrt{d (
\prod_{i=1}^d p_i )^{{1}/{d}}}} \biggr)^N.
\]
Since $\sum_{i=1}^d p_i = 1$, it follows by the AM-GM inequality that
$d (
\prod_{i=1}^d p_i )^{{1}/{d}} < 1$, unless all the entries of
$\pv$ are equal.
Hence, if all the entries of $\mathbf{x}$ are the same and all entries of
$p_i$ are not the same,
the coefficient of $( 1 - \frac{s}{N} )^n$ in the crude
upper bound is
exponential in $N$. If all the $p_i$ are same, the coefficient is of
the order
$N^{{(d - 1)}/{4}}$.

The main fact is that the coefficient of $( 1 - \frac{s}{N}
)^n$ in the upper
bound derived in~(\ref{eqehf2}) is linear in $N$, whereas the
coefficient of $( 1 -
\frac{s}{N} )^n$ in the crude upper bound is almost always
exponential in $N$.
\end{rem*}

\begin{eg}
Consider the generalized Ehrenfest urn model where $N=100$ balls are
distributed in $d=5$
urns. Suppose $s=1$ ball is chosen and each urn is chosen with
probability $p_i=1/5,
i=1,2,\ldots,5$. Using (\ref{eqehf1}) and (\ref{eqehf2}), for a
starting state
$\mathbf{x}=(0,20,0,20,60)$, the bounds on the total variation
distance are
obtained as
%
\begin{equation}\label{eqehf3}
0.25 \biggl( 1 - \frac{1}{100} \biggr)^{n} \leq\| K_\mathbf{x}^n
-\pi
\|_{\mathrm{TV}} \leq
120\biggl ( 1 - \frac{1}{100} \biggr)^{n}.
\end{equation}
For $\varepsilon=0.01$, (\ref{eqehf3}) tells us that $321$ steps are
necessary and $935$
steps are sufficient for the total variation distance to be less then
$0.01$. The
crude upper bound for total variation distance $\frac{1}{2 \sqrt{\pi
(\mathbf{x}
)}} ( 1 -
1/100 )^n =(1.02 \times10^{15}) ( 1 - 1/100 )^n$
would have implied
that $3897$ steps are sufficient for the total variation distance to be
less then $0.01$.
\end{eg}

\section{Discussion} \label{sec4}

We use a probabilistic technique based on a monotone coupling argument
for analyzing all the
examples in this paper. We obtain reasonable upper and lower bounds for
the total variation
distance for any arbitrary starting point of the Markov chain,
significantly broadening previous
results in~\cite{Khare09}. This analysis is very simple to implement,
requiring the knowledge of
a single eigenfunction and its corresponding eigenvalue. In addition,
the analysis does not
require the assumption of reversibility. As an illustration, we provide
the nonreversible
Moran model in Section~\ref{secmoran}. The next goal is to sharpen the
bounds to obtain
matching upper and lower bounds, and to generalize the techniques
developed in this paper for
continuous state spaces.

\begin{appendix}
\section*{Appendix}\label{app}

\newtheorem{lemma}{Lemma}
\begin{lemma}
If the mutation matrix $\uM$ is irreducible, then the transition
density~$K$ in (\ref{mo1})
is irreducible and aperiodic.
\end{lemma}

\begin{pf}
We first show irreducibility. Let $\mathbf{x}\in\mathcal{X}_N^d$ be
arbitrarily chosen. Let $i \neq
j$ be such that $1 \leq i, j \leq d$ and $x_i > 0$. By the
irreducibility of $\uM$, there exists
$n \in\mathbb{N}$ such that $(\uM^n)_{ij} > 0$. As a result, there
exist $i = k_0, k_1, k_2,
\ldots, k_{n-1}, k_n = j$ such that $\prod_{l=0}^{n-1} m_{k_l k_{l+1}}
> 0$. Let
$\mathbf{x}^0 = \mathbf{x}$, and $\mathbf{x}^l = \mathbf{x}^{l-1} +
\ev_{k_l}
- \ev_{k_{l-1}}$ for
$1 \leq l \leq n$. Note
that by construction, $x^l_{k_l} > 0$, which implies $\mathbf{x}^l \in
\mathcal
{X}_N^d$ for every
$1 \leq l \leq n$. Hence,
\begin{eqnarray*}
K^n (\mathbf{x}, \mathbf{x}+ \ev_j - \ev_i)
&=& K^n (\mathbf{x}^0, \mathbf{x}^n)\\
&\geq& \prod_{l=0}^{n-1} K(\mathbf{x}^l, \mathbf{x}^{l+1})\\
&\geq& \prod_{l=0}^{n-1} \frac{x^l_{k_l}}{N} \frac{x^l_{k_l}}{N} m_{k_l
k_{l+1}}\\
&>& 0.
\end{eqnarray*}
We have thus shown that if $\mathbf{x}$ and $\mathbf{y}$ are
neighbors in
$\mathcal
{X}_N^d$, that is, if $\mathbf{y}$
can be obtained from $\mathbf{x}$ by removing an individual in one
species and
adding an individual
in another, then there exists $n \in\mathbb{N}$ such that $K^n
(\mathbf{x},
\mathbf{y}) > 0$. Since any two
elements of $\mathcal{X}_N^d$ are connected by a path such that
successive elements in the
path are neighbors, it follows that $K$ is irreducible.

We now show aperiodicity. Since $\uM$ is irreducible, there exist $i,
j$ such that $1 \leq i
\neq j \leq d$ and $m_{ij} > 0$. If $\mathbf{x}\in\mathcal{X}_N^d$
is such
that $x_i, x_j > 0$,
then
\[
K(\mathbf{x}, \mathbf{x}) \geq\frac{x_j}{N} \frac{x_i}{N} m_{ij} > 0.
\]
Since $K$ is irreducible, and there exists at least one $\mathbf{x}\in
\mathcal{X}_N^d$ such that
$K(\mathbf{x}, \mathbf{x}) > 0$, it follows that $K$ is aperiodic.
\end{pf}
\end{appendix}
%


\printaddresses

\end{document}